\documentclass[11pt, oneside]{article}   	
\usepackage{geometry}                		
\usepackage{graphicx}				
\usepackage{amssymb}

\usepackage{amssymb,amsmath,amsthm}
\newtheorem{thm}{Theorem}
\newtheorem{lem}[thm]{Lemma}

\newtheorem*{ex}{Example}

\newtheorem{conj}{Conjecture}
\newtheorem*{pf}{Proof}

\numberwithin{question}{section}

\title{On Stanley-Reisner rings with linear resolutions}
\author{Ralf Fr\"oberg}
\date{}							

\begin{document}
\maketitle
\begin{abstract}
For a graph $G$, Bayer-Denker-Milutinovi\'c-Rowlands-Sundaram-Xue study in \cite{B-D-M-R-S-X} a new graph complex $\Delta_k^t(G)$,
namely the simplicial complex with facets that are complements to independent sets of size $k$ in $G$. They are interested in topological 
properties such as shellability, vertex decomposability, homotopy type, and homology of these complexes. In this paper we study
more algebraic properties, such as Cohen-Macaulayness, Betti numbers, and linear resolutions of the Stanley-Reisner ring of these complexes and their
Alexander duals.

\end{abstract}

\section{Preliminaries}
A simplicial complex on $[1,n]=\{1,2,\ldots,n\}$ is a collection $\Sigma$ of subsets $F_i$ of $[1,n]$ such that if $F\in\Sigma$ and
$G\subset F$, then $G\in\Sigma$. Maximal faces are called facets. The complex is
called pure if all facets have the same dimension.
If $\Sigma$ is a simplicial complex with vertices $\{x_1,\ldots,x_n\}$, its Stanley-Reisner ideal $I_\Sigma$ is generated by 
the squarefree monomials $x_{i_1}\cdots x_{i_k}$ for all non-faces $\{x_{i_1},\ldots,x_{i_k}\}$
of $\Sigma$. Its Stanley-Reisner ring over a field $k$ is $k[\Sigma]=k[x_1,\ldots,x_n]/I_\Sigma$.

A clique in a graph is a subset $F$ of the vertices of $G$ for which all pairs of vertices in 
$F$ are connected with an edge. The clique complex of $G$ is the simplicial complex which consists of all cliques in $G$.

For a simplicial complex $\Delta$ with vertex set $[1,n]$, its Alexander dual is 
$\{ F\subseteq [1,n];[n]\setminus F\notin\Delta\}$.

A graded algebra $R=k[x_1,\ldots,x_n]/I=S/I$ has a pure resolution if the minimal resolution looks like this.
$$0\leftarrow S/I\leftarrow S\leftarrow S[-a_1]^{b_1}\leftarrow S[-a_2]^{b_2}\leftarrow\cdots\leftarrow S[-a_p]^{b_p}\leftarrow0$$
where $S[-i]$ means that we have shifted $S$ $i$ steps. Thus $S/I$ has a pure resolution if (Tor$^S_i(S/I,k))_j\ne0$ for at most one $j$
for each $i$. If $S/I$ has a pure resolution, the Hilbert
series gives all graded Betti numbers.The Hilbert series of a ring with pure resolution is 
$\frac{1-b_1t^{a_1}+b_2t^{a_2}-\cdots(-1)^pb_pt^{a_p}}{(1-t)^n}$. Now
$S/I$ has a $s$-linear resolution if it has a pure resolution with $a_i=s+i-1$ for each $i$. 
Thus $S/I$ has a $s$-linear resolution if $I$ is generated in degree $s$ and all higher syzygies
are linear. It has a linear resolution if it has a $s$-linear resolution for some $s$. Thus we have

\begin{lem}\label{lin} A ring $k[x_1,\ldots,x_n]/I$ with $s$-linear resolution
$$0\leftarrow S/I\leftarrow S\leftarrow S[-s]^{b_1}\leftarrow S[-(s+1)]^{b_2}\leftarrow\cdots\leftarrow S[-(s+p-1)]^{b_p}\leftarrow0$$
has Hilbert series 
$$\frac{1-b_1t^{s}+b_2t^{s+1}-\cdots+(-1)^pb_pt^{s+p-1}}{(1-t)^n}.$$
\end{lem}

\medskip
We will all the time use the following fundamental theorem by Eagon-Reiner.

\begin{thm}\cite[Theorem 3]{E-R} $k[\Delta]$ has a linear resolution if and only if $k[\Delta^\vee]$ is Cohen-Macaulay.
\end{thm}

For a simplicial complex $\Delta$, let $G(\Delta)$ be its 1-skeleton, i.e., the faces of $\Delta$ of dimension $\le1$. For a graph $G$, 
let $\Delta(G)$ be the largest simplicial complex
with $G$ as 1-skeleton, i.e., $\Delta(G)$ is the clique complex of $G$.  We have that $I_\Sigma$ 
is generated in degree 2 if and only if $\Sigma=\Delta(G)$ for some graph $G$.  The complement graph $\overline G=(V,E')$ to $G=(V,E)$
has edges $E'=\{ (i,j);(i,j)\notin E\}$. The following was proved in \cite{Fr}:

\begin{thm}\cite[Theorem 1]{Fr} Suppose $\Sigma=\Delta(G)$. Then $I_\Sigma$ has a 2-linear resolution if and only if
$G$ is a chordal graph, i.e., every cycle in $G$ of length $>3$ has a chord.
\end{thm}

The following extension was proved by Eagon-Reiner.
\begin{thm}\cite[Proposition 8]{E-R}\label{eq} The following are equivalent for a graph $G$.
\begin{itemize}
\item $\Delta(G)^\vee$ is vertex-decomposable.
\item $\Delta(G)^\vee$ is Cohen-Macaulay over any field $k$.
\item $\Delta(G)^\vee$ is Cohen-Macaulay over some field $k$.
\item $G$ is chordal.
\end{itemize}
\end{thm}

A pure $d$-dimensional simplicial complex is pure shellable if it is pure and one can order its facets $F_1,\ldots,F_m$ such that for 
$i\ge2$ $\langle F_i\rangle\cap(\cup_{j=1}^{i-1}\langle F_j\rangle$, where $\langle H\rangle$ is the simplex generated by $H$, 
is a pure simplicial complex of dimension $d-1$ for all $i$.

\medskip
Let $\sigma$ be a face in the simplicial complex $\Sigma$. The link of $\sigma$ in $\Sigma$ is 
lk$_\Sigma \sigma=\{\tau\in\Sigma;\tau\cap\sigma=\emptyset,\mbox{  and }\tau\cup\sigma\in\Sigma\}$. 
The deletion of $\sigma$ in $\Sigma$ is del$_\Sigma\sigma=\{\tau\in\Sigma;\tau\not\subseteq\tau\}$.

\medskip
A $d$-dimensional simplicial complex $\Sigma$ is vertex decomposable if either $\Sigma$ is a simplex, or there is a vertex $v\in\Sigma$
such that
\begin{itemize}
\item both lk$_\Sigma v$ and del$_\Sigma v$ are vertex decomposable, and
\item del$_\Sigma v$ is pure and $d$-dimensional.
\end{itemize}

It is known that a pure vertex decomposable complex is pure shellable, \cite[Corollary 2.9]{P-B} and a pure shellable complex is Cohen-Macaulay, \cite{Fo}.
Thus one can add $\Delta(G)^\vee$ is pure shellable to the set of equivalences of Theorem \ref{eq} above.

\section{Introduction}
For a graph $G=(V,E)$, the authors of \cite{B-D-M-R-S-X} introduced the simplicial  complex $\Delta_k^t(G)$ whose facets are 
complements to independent sets of size $k$ in $G$. (A set $M$ of vertices is independent if there are no edges between vertices in $M$.)
We have that $\sigma$ is a face of $\Delta_k^t(G)$ if and only if its complement contains an independent set of size $k$. Equivalently, the facets of 
$\Delta_k^t(G)$ are the vertex covers of size $|V|-k$, i.e., sets of vertices that contain at least one endpoint of every edge in $E$.
The Alexander dual $(\Delta_k^t(G))^\vee$  of $\Delta_k^t(G)$ is $\{\sigma\subseteq V; \mbox {no $k$-subset of $\sigma$ is independent in $G$}\}$. 
In particular $(\Delta(G)_2^t)^\vee$ is the clique complex of $G$. The theorem above can thus be formulated, \cite{E-R}:

\begin{thm} The Stanley-Reisner ring of $(\Delta_2^t)(G)^\vee$ has a linear resolution if and only if $G$ is chordal.
\end{thm}

This was extended by Eagon-Reiner in \cite[Proposition 8]{E-R}:

\begin{thm} The graph $G$ is chordal if and only if $\Delta_2^t(G)$ is vertex decomposable and if and only if $\Delta_2^t(G)$ is Cohen-Macaulay.
\end{thm}

We will use another fundamental theorem, this time by Reisner.

\begin{thm}\cite[Theorem 1]{R}\label{reis} Let $\Sigma$ be a simplicial complex on $[1,n]$. Then $k[\Sigma]$ is Cohen-Macaulay 
if and only if $\Sigma$ has Cohen-Macaulay links of vertices and $\tilde H_i(\Sigma;k)=0$ if $i<\dim\Sigma$. In particular, if $\Sigma$
is 1-dimensional (a graph) then $k[\Sigma]$ is Cohen-Macaulay if and only if $G$ is connected.
\end{thm}

The $f$-vector of a simplicial complex $\Sigma$ on $[1,n]$ is $f(\Sigma)=(f_{-1},f_0,\ldots,f_n)$, where $f_i$ is the number of 
$i$-dimensional faces in $\Sigma$. Sometimes we skip the zeros at the end.

\medskip
The three lemmas that follow are certainly well known. The way to determine the Hilbert series of a Stanley-Reisner ring in the first lemma
is easily seen.

\begin{lem}\label{hil}
If $\Sigma$ has $f$-vector $(f_{-1},\ldots,f_d)$, then the Hilbert series of $k[\Sigma]$ is 
$$f_{-1}+\frac{f_0t}{1-t}+\frac{f_1t^2}{(1-t)^2}+\cdots+\frac{f_dt^{d+1}}{(1-t)^{d+1}}.$$
\end{lem}

Also the following lemma is easy.

\begin{lem}\label{dualf}
If the $f$-vector of $\Sigma$ on $[1,n]$  is $(f_{-1},f_0,\ldots,f_n)$, then the $f$-vector of $\Sigma^\vee$ is $(h_{-1},h_0,\ldots,h_n)$, where
$h_i={n\choose i+1}-f_{n-i-2}$.
\end{lem}

Sometimes we may use the following lemma.

\begin{lem}\label{dualgen}
If $\Sigma$ is a simplicial complex on $[1,n]$ with facets $F_1,\ldots,F_k$, then the Stanley-Reisner ideal $I_{\Sigma^\vee}$ is generated by
$\prod_{j\in[1,n]\setminus F_1}x_j,\ldots,\prod_{j\in[1,n]\setminus F_k}x_j$.
\end{lem}

\medskip
In \cite{B-D-M-R-S-X} Bayer-Denker-Milutinovi\'c-Rowlands-Sundaram-Xue studied the complexes $\Delta_k^t(G)$ for from a topological point of view
(shellability, homotopy, homology). Also a lot of specific classes of graphs are considered. In this paper we look at $\Delta_k^t(G)$ and $\Delta_k^t(G)^\vee$ with respect to linear resolution, Betti numbers, and Cohen-Macaulayness of Stanley-Reisner rings. In a second
paper by Bayer-Denker-Milutinovi\'c-Rowlands-Sundaram-Xue \cite{B-D-M-R-S-X2} they study a related graph complex $\Delta_k(G)$ whose facets are
sets of size $k$ such that the induced subgraph is disconnected. Of course $\Delta_2(G)=\Delta_2^t(G)$. In \cite{B-D-M-R-S-X2} a lot of special
classes of graphs are considered with respect to topological properties. This inspired us to consider these classes in our study of algebraic properties.

\section{$k=1$}
We start with an easy lemma which will be used several times. The $i$- dimensional skeleton $\Sigma^{(i)}$ of a simplicial complex 
$\Sigma$ consists of all faces in $\Sigma$ of dimension $\le i$. The lemma and theorem that follows are certainly well known.

\begin{lem}\label{skel}
If $\Sigma$ is an $n$-dimensional simplex, then $\Sigma^{(i)}$ is the Alexander dual to the $\Sigma^{(n-i-2)}$.
The Stanley-Reisner ring of $\Sigma^{(i)}$ is Cohen-Macaulay and has a linear resolution. The Stanley-Reisner ideal of $\Sigma^{(i)}$
is generated by all squarefree monomials of degree $i+2$.
\end{lem}

\begin{pf} $\Sigma^{(i)}$ has $f$-vector $({n+1\choose0},{n+1\choose1},\ldots,{n+1\choose i+1})$. By Lemma \ref{dualf}, 
$(\Sigma^{(i)})^\vee$ has $f$-vector
$({n+1\choose0},{n+1\choose1},\ldots,{n+1\choose n-i-1})$. Thus $(\Sigma^{(i)})^\vee$ is $\Sigma^{(n-i-2)}$.
For the Cohen-Macaulayness we use Theorem \ref{reis}. A skeleton of a simplex has no homology
below the dimension, and links of vertices are Cohen-Macaulay by  induction. Since both $\Sigma^{(i)}$ and its dual are Cohen-Macaulay,
they both have linear resolutions. $\Sigma^{(n-i-2)}$ has facets of dimension $n-i-2$,
so of size $n-i-1$. Thus $(\Sigma^{(n-i-2)})^\vee=\Sigma^{(i)}$ has Stanley-Reisner ideal generated by $n+1-(n-i-1)=i+2$ elements.
\end{pf}

For any graph $G$ with $n$ vertices, we have that $\Delta_1^t(G)$  have all subsets with $n-1$ elements as facets, so it is the $(n-2)$-dimensional
skeleton of an $(n-1)$-dimensional simplex. Its dual is the empty set.

\begin{thm} For any graph with $n$ vertices, $k[\Delta_1^t(G)]=k[x_1,\ldots,x_n]/(x_1x_2\cdots x_n)$ and 
$k[\Delta_1^t(G)^\vee]=k[x_1,\ldots,x_n]/(x_1,\ldots,x_n)$ by Lemma \ref{dualgen} and Lemma \ref{skel}.
They are both Cohen-Macaulay, so they both have linear resolutions. The Betti numbers are
$\beta_{0,0}=1,\beta_{1,n}=1$ and $\beta_{i,i}={n\choose i}$, $i=1,\ldots,n$, respectively.
\end{thm}

\section{$k=2$}

\begin{lem} If $G$ is triangle free, then $\Delta_2^t(G)^\vee=G$.
\end{lem}

\begin{pf} $G$ is the clique complex of $G$.
\end{pf}

It is known which simplicial complexes that have a Stanley-Reisner ring with 2-linear resolution. These are the so called fat forests.
Here is a recursive definition of fat forests, \cite{Fr},  \cite{Fr1}. A simplex $F_1$ is a fat forest. If $F_i$, $i=1,\ldots,k$ are simplices and 
$G_{k-1}=\cup_{I=1}^{k-1}F_i$ is a fat forest, then $G_{k-1}\cup F_k$ is a fat forest if $H=G_{k-1}\cap F_k$ is a simplex, $\dim H\ge-1$.
(If $\dim H=-1$, then $G_{k-1}$ and $F_k$ are disjoint.)

\begin{thm}\cite[Theorem 1]{Fr} $k[\Sigma]$ has a 2-linear resolution if and only if $\Sigma$ is a fat forest (called generalized forest in \cite{Fr}).
\end{thm}

Here is an easy way to determine the Hilbert series of $k[\Sigma]$ if $\Sigma$ is a fat forest.

\begin{thm}\cite[Theorem 1]{Fr1}\label{fat} Let $\Sigma=\Sigma_1\cup\cdots\cup\Sigma_k$ be a fat forest, where $\Sigma_i$ is a simplex of dimension $d_i$,
$i=1,\ldots,k$, and $(\Sigma_1\cup\cdots\cup\Sigma_{j-1})\cap \Sigma_j$ a simplex of dimension $r_i$, $j=2,\ldots,k$. Then the Hilbert series of $k[\Sigma]$
is $\sum_{i=1}^k\frac{1}{(1-t)^{d_i+1}}-\sum_{j=2}^k\frac{1}{(1-t)^{r_i+1}}$.
\end{thm}

\section{$k\ge2$}

We will now look at some special classes of graphs.

\subsection{0-dimensional graphs}
Let $P_n$ consist of $n$ isolated points.

\begin{thm} Both $\Delta_k^t(P_n)$ and its dual have Stanley-Reisner rings which are Cohen-Macaulay and  have linear resolutions.
\end{thm}

\begin{pf}
This follows from Lemma \ref{skel}.
\end{pf}

\begin{ex} $\Delta_2^t(P_n)$ equals the $(n-3)$-dimensional skeleton of an $(n-1)$-dimensional simplex, and  $\Delta_2^t(P_n)^\vee=P^n$
the 0-dimensional skeleton. The $f$-vector of $\Delta_2^t(P_n)$ is $({n\choose 0},{n\choose 1},\ldots,{n\choose n-2})$, and the dual $P^n$
has $f$-vector $(1,n)$. The Hilbert series of $k[\Delta_2^t(P_n)]$ is 
$$1+\frac{nt}{1-t}+\frac{{n\choose2}t^2}{(1-t)^2}+\cdots+\frac{{n\choose n-2}t^{n-2}}{(1-t)^{n-2}}=\frac{1}{(1-t)^n}-\frac{t^n}{(1-t)^n}-\frac{nt^{n-1}}{(1-t)^{n-1}}=$$
$$\frac{1-nt^{n-1}+(n-1)t^{n-1}}{(1-t)^n}.$$ 
The Betti numbers are $\beta_0=\beta_{0,0}=1,\beta_1=\beta_{1,n-1}=n,\beta_2=\beta_{2,n}=n-1$.
The Hilbert series of the dual $k[\Delta_2^t(P_n)]$ is
$$1+\frac{nt}{1-t}=\frac{(1-t)^n+nt(1-t){n-1}}{(1-t)^n}=\frac{1+\sum_{i=1}^{n-1}(-1)^i(n{n-1\choose i}-{n\choose i+1})t^{i+1}}{(1-t)^n}.$$
The Betti numbers are $\beta_0=\beta_{0,0}=1,\beta_i=\beta_{i,i+1}=n{n-1\choose i}-{n\choose i+1}$ for $i=1,\ldots,n-1$,
\end{ex}

\subsection{Trees}
\begin{thm}
Suppose $T$ is a tree with $n$ vertices. Then both $k[\Delta_2^t(T)]$ and $k[\Delta_2^t(T)^\vee]$ are Cohen-Macaulay, and thus have 
linear resolutions. The Betti numbers of $k[\Delta_2^t(T)]$ are $\beta_0=\beta_{0,0}=1,\beta_1=\beta_{1,n-2}=n-1,\beta_2=\beta_{2,n-1}=n-2$, 
and the Betti numbers of $k[\Delta_2^t(T)^\vee]=k[T]$ are
$\beta_{i,i+1}(k[T])=(n-2){n-2\choose i+1}+{n-2\choose i+1}$ for $i=1,\ldots,n-1$.
\end {thm}

\begin{pf} $\Delta_2^t(T)^\vee=T$. That $k[T]$ is Cohen-Macaulay follows from Theorem \ref{reis}. $T$ is a fat forest, so $k[T]$ has a 2-linear resolution.
Thus $\Delta_2^t(T)=(\Delta_2^t(T)^\vee)^\vee$ also has a Cohen-Macaulay Stanley-Reisner ring with linear resolution.
The $f$-vector of $T=\Delta_2^t(T)^\vee$ is $(1,n,n-1)$, so the Hilbert series is 
$$1+\frac{nt}{1-t}+\frac{(n-1)t^2}{(1-t)^2}=\frac{(1-t)^n+nt(1-t)^{n-1}+(n-1)t^2(1-t)^{n-2}}{(1-t)^n}=$$
$$\frac{1+\sum_{i=1}^{n-1}(-1)^{i+1}{n\choose i+1}-n{n-1\choose i}+(n-1){n-2\choose i-1}t^{i+1}}{(1-t)^n}.$$
Thus $\beta_{i,i+1}(k[T])={n\choose i+1}-n{n-1\choose i}+(n-1){n-2\choose i-1}$ for $i=1,\ldots,n-1$.
According to Lemma \ref{dualf} we have $f_i(\Sigma^\vee)={n\choose i+1}-f_{n-i-2}(\Sigma)$, so the $f$-vector of $\Delta_2^t(T)$ is
${n\choose0},{n\choose1},\ldots,{n\choose n-3},{n\choose2}-(n-1),0,0)$, which gives that the Hilbert series is $\frac{1-(n-1)t^{n-2}+(n-2)t^n}{(1-t)^{n-1})}$,
so $\beta_{1,n-2}=n-1,\beta_{2,n-1}=n-2$.
\end{pf}

So the Stanley-Reisner rings of all trees with $n$ vertices have the same Hilbert series, 
in fact even the same Betti numbers. Thus the same is true for
$k[\Delta_2^t(T)]$ for all trees $T$ with $n$ vertices. However, as we will see, two trees with $n$ vertices may have different
Hilbert series of $k[\Delta_3^t(T)]$. 

\subsection{$S_n$}
Let $S_n$ be the graph on $[1,n]$ with edges $\{i,n\}$, $i=1,\ldots,n-1$. 

\begin{thm}
Both $k[\Delta_k^t(S_n)]$ and $k[\Delta_k^t(S_n)^\vee]$ are Cohen-Macaulay, and thus have linear resolutions.
\end{thm}

\begin{pf}
$\Delta_k^t(S_n)$ has facets $\{x_n\}\cup M_i$, where $M_i$ goes through all $(n-k-1)$-subsets of $\{x_1,\ldots,x_{n-1}\}$. 
The Alexander dual has facets $\{x_n\}\cup Q_i$, where $Q_i$ goes through all $(k-1)$-subsets of $\{x_1,\ldots,x_{n-1}\}$.
The Stanley-Reisner ring of $\Delta_k^t(S_n)$ has relations all squarefree monomials in $\{x_1,\ldots,x_{n-1}\}$ of degree $n-k$,
and the dual of $\Delta_3^t(S_n)$ has relations all squarefree monomials in $\{x_1,\ldots,x_{n-1}\}$ of degree $k$.Thus, 
the Stanley-Reisner ring of $\Delta_k^t(S_n)$ is $k[\Sigma][x_n]$, where $\Sigma$ is the $(n-k-2)$-dimensional skeleton
of the simplex on a $(n-2)$-dimensional simplex, and the Stanley-Reisner ring of the dual is $k[\Delta][x_n]$, where $\Delta$ 
is the $(k-2)$-dimensional skeleton of the simplex on a $(n-2)$-dimensional simplex. Stanley-Reisner rings of skeletons of a simplex
are Cohen-Macaulay, since they have no homology below their dimension, and links of vertices are Cohen-Macaulay by induction. 
\end{pf}

\medskip
\begin{ex}
Both $k[\Delta_2^t(S_6)]$ and $k[\Delta_2^t(S_6)^\vee]$ have Betti numbers $\beta_1=\beta_{1,3}=10,\beta_2=\beta_{2,4}=15,
\beta_3=\beta_{3,5}=6$, while $k[\Delta_3^t(L_6)]$ has Betti numbers $\beta_1=\beta_{1,2}=9,
\beta_2=\beta_{2,3}=18,\beta_3=\beta_{3,4}=15,\beta_4=\beta_{4,5}=6,\beta_5=\beta_{5,6}=1$, and its dual $\beta_{0,0}=1,\beta_1=\beta_{1,3}=2,
\beta_2=\beta_{2,6}=1$. Thus $k[\Delta_3^t(S_6)]$ is Cohen-Macaulay and has a linear resolution, while $k[\Delta_3^t(L_6)]$ has a linear resolution,
but is not Cohen-Macaulay.
\end{ex}

\subsection{$K_2\times K_n$}
$K_2\times K_n$ has vertices $\{ x_1,\ldots,x_n\}\cup\{ y_1,\ldots,y_n\}$ and edges 
$\{ (x_i,x_j);1\le i<j\le n\}\cup\{ (y_i,y_j);1\le i<j\le n\}\cup\{(x_i,y_i);1\le i\le n\}$. The maximal independent sets are $\{x_i,y_j;i\ne j\}$.

\begin{thm} $k[\Delta_2^t(K_2\times K_2]$ has a linear resolution with Betti numbers $\beta_{1,2}=\beta_{2,3}=4,\beta_{0,0}=
\beta_{3,4}=1$. Its dual is a complete intersection $k[x_1,x_2,x_3,x_4]/(x_1x_3,x_2,x_4)$ with Betti numbers $\beta_{1,2}=2,
\beta_{0,0}=\beta_{2,4}=1$. If $n>2$, then neither $k[\Delta_2^t(K_2\times K_n)]$ nor $k[(\Delta_2^t(K_2\times K_n))^\vee]$ are 
Cohen-Macaulay, so neither has a linear resolution.
\end{thm}

\begin{pf} First $\Delta_2^t(K_2\times K_2)$ has facets $\{x_1,y_2\},\{x_2,y_1\}$ and its Stanley-Reisner ideal is $(x_1x_2,x_1y_1,x_2y_2,y_1y_2)$.
The dual has facets $\{x_1,x_2\},\{x_2,y_2\},\{y_1,y_2\},\{x_1,y_1\}$ and Stanley-Reisner ideal $(x_1y_2,x_2y_1)$.
$\Delta_2^t(K_2\times K_2)$ is a fat forest, so it has a 2-linear resolution. The nonzero Betti numbers are $\beta_{1,2}=\beta_{2,3}=4,
\beta_{0,0}=\beta_{3,4}=1$. The dual is not a fat forest, but has Cohen-Macaulay Sanley-Reisner ring $k[x_1,x_2,y_1,y_2]/(x_1y_2,x_2y_1)$
with Betti numbers $\beta_{0,0}=\beta_{2,4}=1,\beta_{1,2}=2$. If $n>2$ $k[\Delta_2^t(K_2\times K_n)]$ has minimal relations
$\prod_{i=1}^nx_i,\prod_{i=1}^ny_i,\prod_{i=1}^n(x_iy_i)/x_jy_j,j=1,\ldots,n$. They don't have the same degree, so the resolution
is not linear. The dual has facets $\{x_1,\ldots,x_n\},\{y_1,\ldots,y_n\}$ and $\{x_iy_i\}, i=1,\ldots,n$ and $k[K_2\times K_n]$has
relations of degree 2.. 
This is not a fat forest, so the resolution is not 2-linear. Thus neither $\Delta_2^t(K_2\times K_n)$ nor its dual is Cohen-Macaulay.
\end{pf}

\subsection{$K_{m,n}$}
The complete bipartite graphs $K_{m,n}$ has vertices $\{x_1,\ldots,x_m\}\cup\{y_1,\ldots,y_n\}$, and edges 
$\{x_i,y_j;1\le i\le m,1\le j\le n\}$. We suppose that $m\le n$. $K_{1,n}=S_{n+1}$ is a tree, which is treated, so we suppose that $m\ge2$.

\medskip
We start with a lemma.
\begin{lem}\label{join} Let $\Sigma_1\star\Sigma_2$ be the join of $\Sigma_1$ and $\Sigma_2$. Then $k[\Sigma_1\star\Sigma_2]=
k[\Sigma_1]\otimes_kk[\Sigma_2]$. $k[\Sigma_1\star\Sigma_2]$ is Cohen-Macaulay if and only if both
$k[\Sigma_1]$ and $k[\Sigma_2]$ are Cohen-Macaulay. If both $k[\Sigma_1]$ and $k[\Sigma_2]$ have $a$-linear resolutions, $a>1$,
then $k[\Sigma_1\star\Sigma_2]$ does not have a linear resolution. 
\end{lem}

\begin{pf} The resolution of $k$ over $k[\Sigma_1\star\Sigma_2]$ is the tensor product of the resolution of $k$ over $k[\Sigma_1]$
and over $k[\Sigma_2]$. Thus the "Betti polynomial" 
$B_{\Sigma_1\star\Sigma_2}(t,s)=\sum_{i,j}b_{i,j}(k[\Sigma_1\star\Sigma_2])t^is^{j-i}$ is the product of the Betti polynomials of the factors.
If $B_{\Sigma_1}=1+a_1ts^{a-1}+a_2t^2s^{a-1}+\cdots$ and $B_{\Sigma_2}=1+b_1ts^{a-1}+b_2t^2s^{a-1}+\cdots$, then $B_{\Sigma_1\star\Sigma_2}=
1+ts^{a-1}(a_1+b_1)+t^2((a_2+b_2)s^{a-1}+a_1b_1s^{2a-2})+\cdots$ which is not the Betti polynomial of a ring with linear resolution
since $a-1\ne2a-2$. For Cohen-Macaulayness, see e.g. \cite{Fr2}.
\end{pf}

\begin{thm} If $m\le k\le n$, then $k[\Delta_k^t(K_{m,n})^\vee]$ is Cohen-Macaulay and has a linear resolution.
If $k<m\le n$, then $k[\Delta_k^t(K_{m,n})^\vee]$ is Cohen-Macaulay, but has not a linear resolution.
\end{thm}

\begin{pf}
If $m\le k\le n$, then $k[\Delta_k^t(K_{m,n})^\vee]=k[x_1,\ldots,x_m,y_1,\ldots,y_n]/I$, where $I$ is generated by all
squarefree monomials of degree $k$ in $\{ y_1,\ldots,y_n]$. This is Cohen-Macaulay and has a linear resolution. If $k<m\le n$, 
then $k[\Delta_k^t(K_{m,n})^\vee]$ is the tensor product of $k[x_1,\ldots,x_m]/I$,
$I$ generated by all squarefree monomials of degree $k$, and $k[y_1,\ldots,y_n]/J$, $J$ generated by all squarefree monomials of degree
$k$. Both factors are Cohen-Macaulay with linear resolution, so their tensor product is Cohen-Macaulay but has not a linear resolution.
\end{pf}

\begin{ex}
$k[\Delta_3^t(K_{4,4})]$ is 4-linear with Betti  numbers $\beta_0=\beta_{0,0},\beta_1=\beta_{1,4}=36, \beta_2=\beta_{2,5}=96, \beta_3=\beta_{3,6}=100, \beta_4=\beta_{4,7}=48,
\beta_5=\beta_{5,8}=9.$ The dual has facets $\{ x_i,x_j,y_k,y_l\}$, $1\le i<j\le 4,1\le j<k\le 4$ and ideal 
$(x_ix_jx_k,y_ly_my_p),1\le i<j<k\le m,1\le l<m<p\le n$. The Stanley-Reisner ring is the tensor product of 
$k[x_1,x_2,x_3,x_4]/I$, $I$ generated by all squarefree monomials of degree 3, and $k[y_1,y_2,y_3,y_4]/J$,
$J$ generated by all squarefree monomials of degree 3.
That the resolution is not linear follows from Lemma \ref{join}.
\end{ex}

\begin{ex} Let $G=K_{m,n}$, $2< m\le n$. Then 
$k[\Delta_2^t(G)]$ has a $(m+n-2)$-linear resolution and the Betti numbers are
$\beta_0=\beta_{0,0}=1,\beta_1=\beta_{1,m+n-2}=mn,\beta_2=\beta_{2,m+n-1}=2mn-m-n,
\beta_3=\beta_{3,m+n}=mn-m-n+1$, 
but is not Cohen-Macaulay.
\end{ex}

\begin{pf} $\Delta_2^t(G)^\vee=K_{m,n}$.
The $f$-vector of $K_{m,n}$ is $(1,m+n,mn)$. Thus the $f$-vector of $\Delta_2^t(K_{m,n})$ is $(h_{-1},\ldots,h_{m+n-1})$, where
 $h_i={m+n\choose i+1}$ for $i=-1,0,\ldots,m+n-4$ and $h_{m+n-3}={m+n\choose m+n-2}-mn,h_{m+n-2}=h_{m+n-1}=0$. The
 Hilbert series is $$\frac{1}{(1-t)^{m+n}}-\frac{mnt^{m+n-2}}{(1-t)^{m+n-2}}-\frac{(m+n)t^{m+n-1}}{(1-t)^{m+n-1}}-\frac{t^{m+n}}{(1-t)^{m+n}}=$$
 $$\frac{1-mnt^{m+n-2}(1-t)^2-(m+n)t^{m+n-1}(1-t)-t^{m+n}}{(1-t)^{m+n}}=$$
 $$\frac{1-mnt^{m+n-2}+(2mn-m-n)t^{m+n-1}-(mn-m-n+1)t^{m+n}}{(1-t)^{m+n}}.$$ 
and this gives the claimed Betti numbers.
\end{pf}     

\subsection{$L_n$} 
$L_n$ has vertices $\{x_1,\ldots,x_n\}$ and edges $\{x_1,x_2\},\{x_2,x_3\},\ldots,\{x_{n-1},x_n\}$. Since $L_n$ is a tree, we can assume that
$k>2$ in $\Delta_k^t(L_n)$. For $k>2$ we have.

\begin{thm} $k[\Delta_k^t(L_{2k-1})]$  has a linear resolution, and is Cohen-Macaulay.
\end{thm}
\begin{pf} The Stanley-Reisner ring of the dual is $k[x_1,\ldots,x_{2k-1}]/(x_1x_3\cdots x_{2k-1})$ which is Cohen-Macaulay and has
(a very short) linear resolution.
\end{pf}

\begin{conj}
$k[\Delta_k^t(L_n)]$  has a linear resolution, and is Cohen-Macaulay for all $n\ge2k-1$.
\end{conj}

\subsection{$C_n,W_{n+1}$}
The cycle $C_n$ has vertices $x_1,\ldots,x_n$ and edges $\{x_1,x_2\},\{x_2,x_3\},\ldots,\{x_{n-1},x_n\},\{x_n,x_1\}$.
$\Delta_2^t(C_n)^\vee=C_n$ and $k[C_n]$ is Gorenstein, \cite[Theorem 5.1]{St}. The wheel graph $W_{n+1}$ has vertices
$\{x_1,\ldots,x_{n+1}\}$ and edges  $\{x_1,x_2\},\{x_2,x_3\},\ldots,\{x_{n-1},x_n\},\{x_n,x_1\}$ and $\{x_i,x_{n+1}\}$, $1\le i\le n$.
If the Stanley-Reisner ring of $\Delta_k^t(C_n)$ is $S$, then the Stanley-Reisner ring of $\Delta_k^t(W_{n+1})$ is $S[x_{n+1}]$
and similarly for their duals. Thus $k[\Delta_k^t(C_n)]$ and $k[\Delta_k^t(W_{n+1})]$ have the same Betti numbers. Similarly for
their duals.

\begin{thm} If $G=C_n$, then $k[\Delta_2^t(G)]$ has a linear resolution with Betti numbers 
$\beta_i=\beta_{i,n-2+i-1}={n\choose i+1}-n{n-1\choose i}+n{n-2\choose i-1}$ for $1\le i\le n-1$.
$k[\Delta_2^t(G)^\vee]=k[C_n]$ is Gorenstein with Betti numbers  $\beta_1=\beta_{1,n-2}=n,\beta_2=\beta_{2,n-1}=n,\beta_3=\beta_{3,n}=1$.
The resolution is not linear unless $n=3$, when it is 1-linear.
\end{thm}

\begin{pf} $C_n$ has $f$-vector $(1,n,n)$, so the Hilbert series of 
$k[C_n]$ is 
$$1+\frac{nt}{1-t}+\frac{nt^2}{(1-t)^2}=\frac{(1-t)^n+nt(1-t)^{n-1}+nt^2(1-t)^{n-2}}{(1-t)^n}=$$
$$\frac{1+\sum_{i=1}^{n-1}(-1)^{i+1}({n\choose i+1}-n{n-1\choose i}+n{n-2\choose i-1})t^{i+1}}{(1-t)^n}.$$
The resolution is symmetric since $k[C_n]$ is Gorenstein. This gives that $\beta_0=\beta_{0,0}=\beta_{n-n}=\beta_{n-2,n}=1$
and $\beta_i=\beta_{i,i+1}={n\choose i+1}-n{n-1\choose i}+n{n-2\choose i-1}$ if $1\le i\le n-3$. (A Gorenstein ring with
$\beta_i=\beta_{i,i+1}$ for $1\le i\le p-1$ and $\beta_{0,0}=\beta_{p,p+2}=1$ is sometimes called 2-linear.) The $f$-vector of $C_n^\vee=\Delta_2^t(C_n)$
is $(h_{-1},h_0,\ldots,h_{n-1})$, where $h_i={n\choose i+1}$ if $-1\le i\le n-4$, $h_{n-3}={n\choose n-2}-n,h_{n-2}=h_{n-1}=0.$
Thus the Hilbert series of $k[\Delta_2^t(c_n)]$ is 
$$\frac{1}{(1-t)^n}-\frac{nt^{n-2}}{(1-t)^{n-2}}-\frac{nt^{n-1}}{(1-t)^{n-1}}-\frac{t^n}{(1-t)^n}=
\frac{1-nt^{n-2}(1-t)^2-nt^{n-1}(1-t)-t^n}{(1-t)^n}=$$
$$\frac{1-nt^{n-2}+nt^{n-1}-t^n}{(1-t)^n}.$$
Thus $\beta_1=\beta_{1,n-2}=n,\beta_2=\beta_{2,n-1}=n,\beta_3=\beta_{3,n}=1$.
\end{pf}

\begin{thm}
The Stanley-Reisner ring of the dual of $\Delta_k^t(C_{2k})$ (and for $\Delta_k^t(W_{2k+1}$) has a linear resolution, but is not Cohen-Macaulay.
\end{thm}
\begin{pf}The Stanley-Reisner ideal of the dual of $\Delta_k^t(C_{2k})$ is $(x_1x_3\cdots x_{2k-1},x_2x_4\cdots x_{2k})$.
These are complete intersections, in particular Cohen-Macaulay,
but they do not have linear resolutions.
\end{pf}

\subsection{$C^2_n$}
We denote by $C^2_n$ the graph on $[1,n]$ with edges $\{i,i+1\pmod n\}$ and $\{i,i+2\pmod n\}$.

\begin{thm}
If $n>6$, then $k[\Delta_2^t(C^2_n)]$
does not have a linear resolution and is not Cohen-Macaulay. If $n=6$, it has a linear resolution with Betti numbers
$\beta_{0,0}=\beta_{4,6}=1,\beta_{1,3}=8,\beta_{2,4}=12,\beta_{3,5}=6$, but is not Cohen-Macaulay.
\end{thm}

\begin{pf} Let $n>6$.The dual of $\Delta_2^t(C^2_n)$ is the clique complex of  $C^2_n$.
The clique complex of $C^2_n$ has not Cohen-Macaulay Stanley-Reisner ring, since the dimension
is 2 and $\tilde H_1(C^2_n)\ne0$. The resolution is not linear,
since the clique complex of $C^2_n$ is not a fat forest. If $n=6$ the Stanley-Reisner ring of the clique complex has ideal 
$(x_1x_4,x_2x_5,x_3x_6)$, and is not a fat forest.
Thus it is Cohen-Macaulay, but has not a linear resolution, so the Stanley-Reisner ring of $\Delta_2^t(C^2_n)$ has a linear resolution, 
but is not Cohen-Macaulay. The $f$-vector is (1,6,12,8) so $\Delta_2^t(C^2_n)$ has $f$-vector (1,6,15,12,3), so the Hilbert
 series is $\frac{1-8t^3+12t^4-6t^5+t^6}{(1-t)^6}$, and the Betti numbers are $\beta_0=\beta_{0,0}=1,
 \beta_1=\beta_{1,3}=8,\beta_2=\beta_{2,4}=12,\beta_3=\beta_{3,5}=6,\beta_4=\beta_{4,6}=1.$
\end{pf}

\subsection{$L^2_n$}
We denote by $L^2_n$ the graph on $[1,n]$ with edges $\{i,i+1\}$, $i=1,\ldots,n-1$ and $\{i,i+2\}$, $i=1,\ldots,n-2$.

\begin{thm} $k[\Delta_2^t(L^2_n)]$ is Cohen-Macaulay and has a linear resolution with Betti numbers 
$\beta_{i,i+1}=(-1)^{i+1}({n\choose i+1}-n{n-1\choose i}+(2n-3){n-2\choose i-1}-(n-2){n-3\choose i-2}$ for $1\le i\le n-3$. 
The clique complex of $L^2_n$ has
Stanley-Reisner ring with linear resolution and with Betti numbers $\beta_0=\beta_{0,0}=1,\beta_1=\beta_{1,n-2}=n-2,\beta_2=\beta_{2,n-1}=n-3$
and is Cohen-Macaulay.
\end{thm}

\begin{pf} The clique complex has no homology and the link of a vertex is 2, 3, or 4 isolated points, so the Stanley-Reisner ring is
Cohen-Macaulay. The clique complex $\Sigma$ of $L^2_n$ is a fat tree, so the resolution of $k[\Sigma]$ is 2-linear. 
$\Sigma=\cup_{i=1}^{n-2} \Sigma_i$, where $\dim\Sigma_i=2$ for all $i$, and $\cup_{i=1}^{j-1}\Sigma_i\cap \Sigma_j$ of dimension 1 for all $j$.
According to Theorem \ref{fat} the Hilbert series is
$$\frac{n-2}{(1-t)^3}-\frac{n-3}{(1-t)^2}=\frac{(1-t)^{n-3}(1+(n-3)t)}{(1-t)^n}.$$
The coefficient of $t^{i+1}$ in the numerator is $(-1)^i({n-3\choose i+1}-(n-3){n-3\choose i}$. We get $\beta_0=\beta_{0,0}=1,
\beta_i=\beta_{i,i+1}=(n-3){n-3\choose i}-{n-3\choose i+1}$, $i=1,\ldots,n-3$.
The $f$-vector of $L^2_n$ is $(1,n,2n-3,n-2)$. This gives the $f$-vector of $\Delta_2^t(L^2_n)$ by Lemma \ref{dualf} and
 the Hilbert series 
 $$\frac{1-(n-2)(1-t)^3t^{n-3}-(2n-3)(1-t)^2t^{n-2}-n(1-t)t^{n-1}-t^n}{(1-t)^n}$$
 of $k[\Delta_2^t(L^2_n)]$ by Lemma \ref{hil}, which gives the Betti numbers since the resolution is linear.
The result is that the resolution is $(n-2)$-linear resolution with $\beta_{1,n-2}=n-2$, $\beta_{2,n-1}=n-3$.
\end{pf}

\begin{conj}
Both $\Delta_k^t(L^2_n)$ and its dual have Cohen-Macaulay Stanley-Reisner rings with linear resolutions.
\end{conj}

\begin{ex}
$k[\Delta_3^t(L^2_8)]$ has Betti numbers $\beta_0=\beta_{0,0}=1,\beta_1=\beta_{1,2}=6,\beta_2=\beta_{2,3}=8,\beta_3=\beta_{3,4}=3$,
and $k[\Delta_3^t(L^2_8)^\vee]$ has Betti numbers $\beta_0=\beta_{0,0}=1,\beta_1=\beta_{1,3}=4,\beta_2=\beta_{2,4}=3$,
\end{ex}

\begin{thm} $k[\Delta_k^t(C^2_{3k})]$ and $k[\Delta_k^t(L^2_{3k-2})]$ have linear resolutions. $k[\Delta_k^t(C^2_{3k})]$ is not
Cohen-Macaulay, but $k[\Delta_k^t(L^2_{3k-2})]$ is.
\end{thm}

\begin{pf} The Stanley-Reisner ring of their duals are
$$k[x_1,\ldots,x_{3k}]/(x_1x_4\cdots x_{3k-2},x_2x_5\cdots x_{3k-1},x_3x_6\cdots x_{3k})\mbox{ and }$$
$$k[x_1,\ldots,x_{3k-2}]/(x_1x_4\cdots x_{3k-2}),$$
respectively which are Cohen-Macaulay, $k[x_1,\ldots,x_{3k}]/(x_1x_4\cdots x_{3k-2},x_2x_5\cdots x_{3k-1},x_3x_6\cdots x_{3k})$
has not a linear resolution, but $k[x_1,\ldots,x_{2k-2}]/(x_1x_4\cdots x_{3k-2})$ do.
\end{pf}

\begin{ex} $\Delta_4^t(C_8)$ has Betti numbers $\beta_0=\beta_{0,0}=1,\beta_1=\beta_{1,2}=16,\beta_2=\beta_{2,3}=48,\beta_3=\beta_{3,4}=68,
\beta_4=\beta_{4,5}=56,\beta_5=\beta_{5,6}=28,\beta_6=\beta_{6,7}=8,\beta_7=\beta_{7,8}=1$.

\medskip\noindent
$\Delta_4(L_7)$ and $\Delta_4^t(L^2_{10})$ has Betti numbers $\beta_0=\beta_{0,0}=1,\beta_1=\beta_{1,1}=4,\beta_2=\beta_{2,2}=6,\beta_3=\beta_{3,3}=4,
\beta_4=\beta_{4,4}=1$.

\medskip\noindent
$\Delta_3^t(C^2_9)$ has Betti numbers $\beta_0=\beta_{0,0}=1,\beta_1=\beta_{1,3}=27,\beta_2=\beta_{2,4}=81,\beta_3=\beta_{3,5}=108,
\beta_4=\beta_{4,6}=81,\beta_5=\beta_{5,7}=36,\beta_6=\beta_{6,8}=9,\beta_7=\beta_{7,9}=1$.

\end{ex}

\subsection{$G_{m,n}$}
The grid $G_{m,n}$ has $mn$ vertices $\{(i,j);1\le i\le m,1\le j\le n\}$ and edges $\{(i,j),(i+1,j)\},1\le i\le m-1,1\le j\le n$ and 
$\{(i,j),(i,j+1)\},1\le i\le m,1\le j\le n-1$. 

\begin{thm} $k[\Delta_2^t(G_{m,n})]$ has a linear resolution with Betti numbers $\beta_0=\beta_{0,0},\beta_1=\beta_{1,mn-2}=2mn-m-n,
\beta_2=\beta_{2,mn-1}=3mn-2m-2n, \beta_3=\beta_{3,mn}=mn-m-n-1$, but is not Cohen-Macaulay. $k[G_{m,n}]$
is Cohen-Macaulay, but has not a linear resolution.
\end{thm}

\begin{pf}
The dual of $\Delta_2^t(G_{m,n})$ is $G_{m,n}$. $G_{m,n}$ has $f$-vector $(1,mn,2mn-m-n)$. 
$k[G_{m,n}]$ is Cohen-Macaulay since $G_{m,n}$ is 1-dimensional and
$\tilde H_0(G_{m,n})=0$ and the link of a vertex is 2, 3, or 4 isolated points. $G_{m,n}$ is not a fat forest, so the resolution of
$k[G_{m,n}]$ is not linear. $k[\Delta_2^t(G_{m,n})]$ has Hilbert series 
$$\frac{1-(2mn-m-n)t^{mn-2}(1-t)^2-mnt^{mn-1}(1-t)-t^{mn}}{(1-t)^{mn}}{(1-t)^{mn}}=$$
$$\frac{1-(2mn-m-n)t^{mn-2}+(3mn-2m-2n)t^{mn-1}+(mn-m-n-1)t^{mn}}{(1-t)^{mn}}.$$
This gives the claimed Betti numbers.
\end{pf}

\noindent
There is no conflict of interest. All data is created by the author.

\end{document}